\documentclass[11pt]{article}

\usepackage[english]{babel}
\usepackage{amsmath,amssymb,enumerate}

\newtheorem{thm}{Theorem}[section]
\newtheorem{cor}[thm]{Corollary}

\newcommand{\N}{{\mathbb N}}

\DeclareMathOperator{\vol}{vol}
\DeclareMathOperator{\iso}{Isom}

\title{Length spectrum of geodesic loops in manifolds of nonpositive curvature}
\author{Muetzel, Bjoern\thanks{E-mail address : bjorn.mutzel@gmail.com}\\
 \\
\small Department of Mathematics, Ecole Polytechnique F\'{e}d\'{e}rale de Lausanne, Station 8, \\
\small CH-1015 Lausanne, Switzerland \\[-0.8ex]
}

\begin{document}

\maketitle

\begin{abstract}
In section 1 we reformulate a theorem of Blichfeldt in the framework of manifolds of nonpositive curvature.
As a result we obtain a lower bound on the number of homotopically distinct geodesic loops
emanating from a common point $q$ whose length is smaller than a fixed constant. This bound depends only
on the volume growth of balls in the universal covering and the volume of the manifold itself.
We compare the result with known results about the asymptotic growth rate of closed geodesics and
loops in section 2.
\end{abstract}

\section{Number of geodesic loops of length smaller than a fixed constant}
We assume in the following that a manifold $M$ is a connected complete Riemannian manifold of nonpositive (sectional) curvature. By covering theory, the universal covering space of $M$, $K$ is a simply connected space of nonpositive curvature. There exists a group $G$ of Deck transformations with the following properties :
\[
   M \simeq {K \mod G}, \text{ \ \ \ } G \subset \iso(K)  \text{ \ and \ } G \simeq \pi_1(M).
\]
A geodesic loop $\delta \subset M$ with base point $q$ defines a homotopy class $[\delta] \in \pi_1(M,q)$. We denote by $l([\delta])$, the \textit{length of $[\delta]$}, which we define as the length of the shortest geodesic loop in the homotopy class of $\delta$ :
\[
  l([\delta])= \min\limits_{\delta' \in [\delta]} l(\delta') .
\]
It is a well known fact that - as the curvature of $M$ is nonpositive - there exists a unique shortest loop in each homotopy class.
We call a non-trivial geodesic loop $\delta$ \textit{primitive}, if
\[
[\delta] \neq [\alpha]^n, \text{ for }  [\alpha] \in \pi_1(M,q), \text{ with } n \in \N \backslash \{1\}.
\]
We say that a geodesic loop $\eta$ is not primitive, if it is a \textit{multiple} of another geodesic loop. We define \textit{primitive} and \textit{multiple} analogously for closed geodesics and their free homotopy classes.

We denote by $v(t)$ the number of primitive (closed) geodesics in $M$ of length smaller or equal to $t$ :
\[
  v(t)= \#\{\gamma \subset M \mid \gamma \text{ primitive geodesic, such that } l(\gamma) \leq t\}.
\]
For $q \in M$, we denote furthermore
\[
P_{q}(t)= \#\{\delta \subset M \mid \delta \text{ geodesic loop with base point } q, \text{ such that } l([\delta]) \leq t \}.
\]

In \cite{bli}, Blichfeldt gives a lower bound on the number of distinct lattice vectors of length smaller than a fixed constant. This bound depends only on the volume of the torus corresponding to the lattice. The theorem relies only on the pigeonhole principle and therefore it can be easily generalized. The following theorem and the subsequent corollaries are a direct generalization of Blichfeldt's theorem in the framework of manifolds of nonpositive curvature.

\begin{thm} Let $M$ be a manifold of nonpositive curvature with $\vol(M) < \infty$. Let $K$ be its universal covering and $p:K\rightarrow M$ a covering map. If $C \subset K$ is a convex set, such that for a $m \in \N$
\[
\vol(C) > m \vol(M),
\]
then there exist pairwise distinct points $x_1,...,x_m \in C$ that map to the same point $q$ in $M$.  The geodesic arcs $(\gamma_{x_i,x_1}) \subset C$ map to $m-1$ homotopically distinct loops $p(\gamma_{x_i,x_1})$ with base point $q$.
\label{thm:Blich_npc}
\end{thm}

Let $B_r(x) \subset K$ be a ball of radius $r$ with center $x$. As $K$ has nonpositive curvature, the ball is convex. Hence we obtain the following corollary from \textbf{Theorem \ref{thm:Blich_npc}}.

\begin{cor} Let $B_r(x)$ be a ball in the universal covering $K$ of $M$, such that for a $m \in \N$, $\vol(B_r(x)) > m \vol(M)$. Then there exist a point $q \in M$ and non-trivial geodesic loops $\delta_2,...,\delta_{m} \subset M $, with common base point $q$ such that
\[
l({\delta_i}) \leq 2r \text{ \ \ for all \ \ } i \in \{2,..,m\}.
\]
Furthermore these geodesic loops are all in different homotopy classes.
\label{thm:Blich_cor}
\end{cor}

In the free homotopy class of each geodesic loop $\delta_i$ of the corollary, there is a unique closed geodesic $\gamma_i$. We note that we can not conclude from the theorem that the $\left(\gamma_i \right)_{i=2,...,m}$ are pairwise distinct. We can also not exclude that such a geodesic $\gamma_i$ is the multiple of another geodesic.\\
As a direct consequence we obtain :
\begin{cor} Let $M$ be a manifold of nonpositive curvature, such that $M$ has finite volume $\vol(M)$. Then there exists a $x \in K$ and a $q = p(x) \in M$, such that
\[
     P_q(t) \geq  \frac{\vol(B_{\frac{t}{2}}(x))}{\vol(M)} -2
\]
\label{thm:Blich_asymptotic}
\end{cor}

\textbf{proof of Theorem~\ref{thm:Blich_npc}} Lift $M$ into the simply connected universal covering $K$ and let $G \subset \iso(K)$ be the group of Deck transformations. Let $F(M)$ be a fundamental domain of the group action of $G$. For $g  \in G$ define
\[
C_g = C \cap g(F(M)).
\]
We have that $C$ is the disjoint union of the $(C_g)_{g \in G}$ and $g^{-1}(C_g) \subset F(M)$.\\
To prove the theorem, we have to show that there exists a point $x \in F(M)$ and pairwise distinct
\[
(g_i)_{i \in \{1,..,m\}} \subset G, \text{ \ \ such that \ \ } x_i=g_i(x) \in C_{g_i}.
\]
It follows then from the convexity of $C$ that the geodesic arc $\gamma_{x_i,x_1}$ connecting $x_i$ and $x_1$ is in $C$. As the $(x_i)_{i=2,..m}$ are pairwise distinct, the $(\gamma_{x_i,x_1})_{i=2,..m}$ map to homotopically pairwise distinct geodesic loops $(\delta_i)_{i=2,..m}$ under the covering map $p$. Their common base point is $p(x_1)$.\\

We show that there exists a $x \in F(M)$ that is covered at least $m$ times by pairwise distinct $g_i^{-1}(C_{g_i})_{i \in \{1,..,m\}}$. We prove the result by contradiction. If every point in $F(M)$ was covered only $m$ times by the disjoint $(g^{-1}(C_g))_{g \in G}$, then we would obtain :
\[
   \vol(C) = \sum\limits_{g \in G} {\vol(C_g)} = \sum\limits_{g \in G} {\vol(g^{-1}(C_g))}
   \leq \int\limits_{F(M)} {m \, d\mu } \leq m \vol(M).
\]

But $\vol(C) > m \vol(M)$, hence we obtain a contradiction.\ \ $\square$\\
\\
We will compare this result with known results about the length spectrum of closed geodesics and geodesic loops in the following section.

\section{Asymptotic growth rate of closed geodesics and loops}

Let $M$ be a manifold of nonpositive curvature and $K$ its universal covering. Let $p : K \rightarrow M$ be the universal covering map. For $q \in M$, we chose a $x \in p^{-1}(q) \subset K$. The \textit{volume entropy} of $M$ is defined as
\[
h_{vol}(M) = \mathop{\lim} \limits_{t \to \infty} \frac{\log(\vol(B_t(x)))}{t}
\]
If $M$ is compact, the limit exists and does not depend on the point $q \in M$  (see \cite{ma}).
If $M$ is a $n$ dimensional manifold, such that the curvature of $M$ is bounded between $-{K_1}^2$ and $-{K_2}^2$, where $0\leq K_1 \leq K_2$, then
\begin{equation*}
(n-1)\cdot K_1   \leq  h_{vol}(M) \leq (n-1)\cdot K_2
\label{eq:vol_K}
\end{equation*}

\textbf{Corollary \ref{thm:Blich_asymptotic}} is related to the following result stated in \cite{kn}.

\begin{thm} Let $M$ be a compact manifold of nonpositive curvature, such that $M$ contains a geodesic, which is not the boundary of a flat half-plane. Then we have for all $q \in M$ that
\[
   \lim\limits_{t \to \infty} \frac{\log(v(t))}{t} = \lim\limits_{t \to \infty} \frac{\log(P_q(t))}{t} = h_{vol}(M).
\]
\label{thm:knieper}
\end{thm}

In comparison, \textbf{Corollary \ref{thm:Blich_asymptotic}} implies the following. If $M$ is a compact manifold then it follows that
\[
   \lim\limits_{t \to \infty} \frac{\log(P_q(t))}{t} \geq \frac{1}{2} h_{vol}(M).
\]
Therefore we underestimate $\lim\limits_{t \to \infty} \frac{\log(P_q(t))}{t}$, if we use the corollary.\\
\textbf{Theorem \ref{thm:knieper}} is additionally stronger in as far, as it makes also a statement about the number of prime geodesics and not only about geodesic loops emanating from a fixed base point. However it treats the asymptotic limit and applies only to compact manifolds, whereas \textbf{Corollary \ref{thm:Blich_asymptotic}} gives a lower bound for the number of homotopically different geodesic loops of length smaller than a certain constant and does not rely on the compactness of $M$.\\
\textbf{Theorem \ref{thm:knieper}} has been generalized to Hadamard manifolds without the necessity of the assumption of compactness or even finite volume of the manifold in \cite{li}.\\

\end{document}